\documentclass[12pt]{article}
\usepackage{amssymb,amsmath}
\usepackage[letterpaper,hmargin=1in,vmargin=1.25in]{geometry}

\newtheorem{thm}{Theorem}

\newtheorem{lem}[thm]{Lemma}

\newtheorem{prop}[thm]{Proposition}

\newtheorem{alg}{Algorithm}

\newcommand{\qed}{\hfill\rule{2mm}{2mm}}

\def\tm#1{\left(\text{mod }#1\right)}

\title{Small
primitive roots and malleability of RSA moduli}

\author
{Luis Dieulefait$^\dag$, Jorge Jim\'enez Urroz$^*$ \\
{\small $^\dag$ Dept. \'Algebra. Universitat de Barcelona}\\
{\small $^*$Dept. Matem\`atica Aplicada IV. Universitat Polit\`ecnica de Catalunya}\\
{\small Campus Nord, c/Jordi Girona, 1-3, 08034 Barcelona}\\
 {\small  e-mail: {ldieulefait@ub.es, jjimenez@ma4.upc.edu }}}
\date{}

\begin{document}


\maketitle

\begin{abstract}
In their paper \cite{pv}, P. Paillier and J. Villar make a
conjecture about the malleability of an RSA modulus. In this paper
we present an explicit algorithm refuting the conjecture. Concretely
we can factorize an RSA modulus $n$ using very little  information
on the factorization of a concrete $n'$ coprime to $n$. However, we
believe the conjecture might be true, when imposing some extra
conditions on the auxiliary $n'$ allowed to be used.  In particular,
the paper shows how subtle the notion of malleability is.

\


\end{abstract}

\section{Introduction}

The existence of a tradeoff between one-wayness and chosen
ciphertext security dates back to the eighties when, for example, it
was observed in \cite{R, W, GMR}. In some sense, one cannot achieve
one-way encryption with a level of security equivalent to solve
certain difficult problem, at the same time as the cryptosystem
being IND-CCA secure respect to it. This so called paradox has been
attempted to be formally proved many times, by a number of authors,
since first observed. However no one succeeded  until very recently,
when Pailler and Villar (cf. \cite{pv}) clarified the question for
the case of factoring-based cryptosystems. In particular, they give
precise conditions for certain security incompatibilities to exist.
More precisely, they reformulate the paradox in terms of key
preserving black-box reductions and prove that if factoring can be
reduced in the standard model to breaking one-wayness of the
cryptosystem then
 it is impossible to achieve chosen-cyphertext security. As the
 authors mention in their paper (cf. \cite{pv}), combining this result with
 the security proofs contained
 in \cite{BR,B} gives a very interesting separation result between the
 Random Oracle model and the standard model. \\

Moreover, assuming an extra hypothesis, which they call
``non-malleability" of the key generator, they are able to extend
the result from key preserving black box reductions to the case of
arbitrary black box reductions. \\
 Hence, as the authors themselves stress in \cite{pv}, it is very
important to study non-malleability of key generators. In fact,
they conjecture that most instance generators are non-malleable,
but no arguments are given to support this belief. The goal of
this note is to shed some light on this open question.\\

Actually, the notion of non-malleability captures a very basic fact
in arithmetic: intuitively, one tends to believe that the problem of
factoring a given number $n$ (an RSA modulus) is not made easier if
we know how to factor other numbers $n'$ relatively prime to $n$.
The random behavior of prime numbers, observed many times in
 the literature, suggests that if the numbers $n'$ are randomly
 selected their factorization is useless for the problem of
 factoring $n$. However this might not be so relevant to malleability because
 we have the freedom to select cleverly the additional numbers
 $n'$.\\

Indeed, the result contained in this note goes against the
non-malleability intuition, thus showing how subtle this notion is.
Concretely, for any number $n$ we are able to prove the existence of
a polynomial time reduction algorithm from factoring
 $n$ to factoring certain explicit numbers $n'$, all relatively
 prime to $n$. In other words, we show that factoring is, in this generality, a malleable problem. \\
 Let us stress that this might be compatible with the conjecture
 of \cite{pv} mentioned above because imposing extra conditions on the
 numbers $n'$ may result in transforming the problem in a
 non-malleable one. In fact, it is our belief that malleability is
 a notion that depends strongly of these kind of extra conditions,
 and hence requires further research.\\





\section{The algorithm}

Given an RSA modulus $n=pq$, we want to find $n'$ such that
factoring $n'$, with the help of an oracle, will allow  us in
finding the factorization of $n$. In fact we will only need very
partial information about the factorization of $n'$ in order to get
the complete factorization of $n$. From now on, and without loose of
generality, we will make the assumption that  $p<q$.
\subsection{A particular case}\label{part}
By construction, (which will be clear in a moment), it turns out
that the particular case  in which $n$ is such that
$2^{p-1}\not\equiv 1\tm{q}$ or $2^{q-1}\not\equiv 1\tm {p}$ is
somehow simpler and we will dedicate this section to it. However,
the whole idea of the method will arise in this case and so the
general one, considered in the next section, will be very similar.
We consider $n'=2^n+1$. Observe that  an efficient encoding of $n'$
of  size comparable to $n$ is available since all these numbers in
binary form have a $1$ at the beginning and end, and the rest are
precisely $n-1$ zeros. Let us assume the existence of an oracle
$\cal O$ which, on input $n'$, returns the residue class modulo $n$
of three prime factors $r|n'$. In fact, the only thing we need is
the residue class of just one factor of $n'$ modulo $n$ different
from $1$ and $3$ so, if convenient, one  can admit an oracle
answering any  set $S\subset\{r\tm{n}\,:\,r\text{ prime }, r|n'\}$,
$S\not\subseteq\{1,3\}$ and polynomial size. We now present an
algorithm which on the input and RSA modulus $n$ in the conditions
of this section, outputs a nontrivial factor of $n$.

\begin{alg}\label{det}

\

\begin{itemize}
\item Send $n'=2^n+1$ in binary form to $\cal O$.
\item Take $r\in S$, $r\ne 1,3$, and compute $d=(r-1,n)$.
\end{itemize}
\end{alg}
\begin{thm}\label{male}
Let $n=pq$ be and RSA modulus such that either $2^{p-1}\not\equiv
1\tm q$ or $2^{q-1}\not\equiv 1\tm p$. Then the  number $d$ given by
the previous algorithm, in polynomial time in $\log n$, is a prime
divisor of $n$.
\end{thm}
{\it Proof:} The first thing we have to prove is that there exists a
set $S$ satisfying the conditions of the algorithm. In order to do
so we have to prove that at least one prime factor of $n'$ is not
$1$ or $3$ modulo $n$.  Suppose $r$ is a prime factor of $n'$. Then,
$2^{2n}\equiv 1\tm r$ and so, either $r=3$ which always divides
$n'$, or  the order of $2$ in $\mathbb F_r^*$ is
ord$_r(2)=p,q,2p,2q,pq$ or $2pq$.  In this case we just have to
recall that the order of any element must divide the order of the
group to conclude that either $p|(r-1)$, $q|(r-1)$ or $n|(r-1)$.
Note, on the other hand that $9$ never divides $n'$ since
$n\equiv\pm1\tm 6$ and so $2^n\equiv 2$ or $5$ modulo $9$. Hence, If
$n|(r-1)$ for any $r|n'/3$,  then each factor of $n'/3$ is $1$
modulo $n$ and so $n'/3\equiv1\tm n$ which is the same as saying
$2^{n-1}\equiv 1\tm n$. This is impossible since in particular
$2^{n-1}\equiv2^{p-1}\tm q$ and $2^{n-1}\equiv2^{q-1}\tm p$. Hence
there exists $r_0|n'$ such that $r_0\not\equiv 1\tm n$. Observe also
that any such factor verifies $r_0\equiv 1\tm p$ or $r_0\equiv 1\tm
q$ and,  in particular, $r_0\not\equiv 3\tm n$.

\qed

The previous algorithm would work, in particular, for any modulus
$n=pq$ such that $(p-1,q-1)=D$ is small, for example $D<\log_2(n)$.
Indeed,  if $2^{p-1}\equiv 1\tm  q$ and $2^{q-1}\equiv 1\tm  p$,
then $2^D\equiv 1\tm n$ which is impossible for $D<\log_2(n)$. This
fact leads to the interesting observation that even the probability
that $D>\log_2(n)$ tends to zero with $n$. This is the content of
the following proposition
\begin{prop}\label{bdh} For any positive $z$ we have
$$
\sum_{{z\le p,q<2z}\atop{(p-1,q-1)>\log z}}1\le \left(\frac{z}{\log
z}\right)^2\frac{(\log\log z)^2}{\log z},
$$
where the sum runs over the prime numbers in the interval.
\end{prop}
{\bf Remark:} Before proving the proposition, let us observe that we
just have to use the Prime Number Theorem to obtain $\sum_{z\le
p,q<2z}1\sim (z/\log z)^2$ and hence, the probability of finding a
pair of primes in the interval $[z,2z]$ which do not satisfy the
conditions in Theorem \ref{male} tends to zero faster than
$(\log\log z)^2/\log z$. Also  note that even if $(p-1,q-1)$ would be big,  we
still would need $2$ to have order $D$ modulo $p$ and modulo $q$ which
one expects to be false for many pairs of primes by Artin's
conjecture, (cf. \cite{m}).

\

\noindent {\it Proof of Proposition} \ref{bdh}. Given a positive $z$
big enough,  let
$$
\pi(d;z)=\sum_{{p\equiv 1\tm{d}}\atop{z\le p< 2z}}1.
$$
Then, the
number of pairs of primes $z\le p,q<2z$ such that $(p-1,q-1)=d>\log
z$ is bounded above by
$$
\sum_{\log z<d<z}\sum_{{p,q\equiv 1\tm{d}}\atop{z\le q<p< 2z}}1<
\sum_{\log
z<d<z^\alpha}\pi(d;z)^2+\sum_{z^\alpha<d<z}\pi(d;z)^2=S_1+S_2,
$$
for any $0<\alpha<1$. For the second term we get trivially the bound
$S_2<4z^{3-2\alpha}$. To estimate $S_1$ let us first introduce the following
useful notation. We will write $E(d;z)=\pi(d;z)-z/(\varphi(d)\log z)$,  as the error in the approximation
of the number of primes in the congruence $1$ modulo $n$ by the total number of
primes divided by the number of congruences. Then,
\begin{eqnarray*}
&&S_1=\sum_{\log z<d<z^\alpha}\left(\frac{z}{\varphi(d)\log
z}+E(d,z)\right)^2 = \\
&&\left(\frac{z}{\log z}\right)^2\sum_{{\log
z<d<2z}}\frac{1}{\varphi(d)^2}+\sum_{\log z<d<z^\alpha}(E(d,z))^2+
2\sum_{\log z<d<z^\alpha}\frac{z}{\varphi(d)\log
z}E(d,z).
\end{eqnarray*}
We can use now  Cauchy-Schwartz inequality to get, for the last sum above
\begin{equation}\label{cs}
\sum_{\log z<d<z^\alpha}\frac{z}{\varphi(d)\log
z}E(d,z)\le
\left(\sum_{\log z<d<z^\alpha}\left(\frac{z}{\varphi(d)\log
z}\right)^2\right)^{1/2}\left(\sum_{\log z<d<z^\alpha}(E(d,z))^2\right)^{1/2}.
\end{equation}
We are in the correct position to  use the
Barban-Davenport-Halberstam Theorem for primes in artihmetic
progressions, (cf. page 421, \cite{ik}), which we now include for
convenience.
\begin{thm} (Barban-Davenport-Halberstam) We have
$$
\sum_{d\le z^{1-\varepsilon}}\left(E(d;z)\right)^2\ll
z^2/(\log z)^A,
$$
for any $A>0$, and  $\varepsilon>0$, where the implied constant only depends on $A$ and $\varepsilon$.
\end{thm}
Substituting the above inequality in $S_1$, putting  $A=4 and
\varepsilon =1/4$, and using (\ref{cs}) we get for some constant $C$
$$
S_1\le
\left(\frac{z}{\log z}\right)^2\sum_{{\log
z<d<2z}}\frac{1}{\varphi(d)^2}+ \frac{Cz^2}{(\log z)^4}+ \frac{Cz}{(\log z)^3}
\left(\sum_{{\log
z<d<2z}}\frac{1}{\varphi(d)^2}\right)^{1/2}.
$$
To finish the proof of the Proposition we just have to note that
$$
\varphi(d)={d\prod_{p|d}(1-1/p)}>d\prod_{p<d}(1-1/p)>Cd/\log d,
$$
by  Mertens Theorem (cf. p.34, \cite{ik}) and so
$$
\sum_{{\log
z<d<2z}}\frac{1}{\varphi(d)^2}\le C\sum_{\log z<d}\left(\frac{\log d}{d}\right)^2\le C_1\frac{(\log\log z)^2}{\log z},
$$
for some constants $C,C_1$. The result follows.

\subsection{ The general case}\label{}
 For a few pairs of primes, it could happen that the order of $2$ in
$\mathbb F_q^*$ and $\mathbb F_p^*$ was a divisor of $D$ and, in
that case, $2^n$ is indeed $2$ modulo $n$ which could make Algorithm
\ref{det} fail. To avoid this problem, instead of $2$, we will
choose a primitive root of $\mathbb F_q^*$, $g$, to build our test
number $n'=g^n+1$. It is very easy to see that the number of
primitive roots of $\mathbb F_q^*$ is $\phi(q-1)$, hence, the
probability for an integer $m$ to be a primitive root verifies 
$$
\frac{\varphi(q-1)}{q-1}=\prod_{p|(q-1)}\left(1-\frac{1}p\right)>\prod_{p<q}\left(1-\frac{1}p\right)\sim\frac{e^{-\gamma}}{\log
{q}},
$$
again by Mertens theorem. In other words, a set of size $C\log q$ of
integers contains a primitive root modulo $q$  with probability as
close to one as we want, making the constant $C$ big enough. To see
this, note that the probability for a random set of this size to
contain no primitive roots would be $(1-1/(e^{\gamma}\log
{q}))^{C\log q}\sim e^{-C/e^{\gamma}}$. In this sense Bach, in \cite{ba}, made a much more accurate heuristic argument to claim that the least primitive root modulo $p$, which we will call $g(p)$ should verify $g(p)\le e^{\gamma}\log p(\log\log p)^2(1+\varepsilon)$ for almost all $p$. Although this fact is not
yet proved, there are conditional results which certify the truth of
the statement. In particular we will mention the following result of
V. Shoup in \cite{Sh} proved under the Grand Riemann Hypothesis, GRH
from now on.
\begin{thm}\label{sho}(Shoup) Let $p$ be a prime and denote $g(p)$ as the least positive integer
which is a generator of  $\mathbb F_p^*$. Then, if GRH is true,
$g(p)=O((\log p)^6)$.
\end{thm}
Observe that, although far from the expected result, $g(p)=O((\log
p)^6)$ is still of polynomial size and, hence, good enough for our
purposes. It is worth mentioning that Heath-Brown was able to prove in \cite{hb} that among $2,3, 5$ there is a primitive root for infinitely many primes $p$. Let us now describe the algorithm.

\

For convenience we will call $c\in\mathbb \{0,1\}^{n+2}$ the binary
encoding of $2^n+1$. We will take advantage of the fact that the
$m$-ary representation of the numbers $m^n+1$ is always $c$,
independent of $m$. Let $n_m'=m^n+1$ and consider the function
$\omega(n)$ counting the number of distinct prime factors of $n$.
Assume the existence of an oracle $\cal O$ which, on input $(c,m)$,
returns a set of residue classes $S$ of size $|S|=\omega(m)+2$ when
such a set $S$ exists, and otherwise returns $\perp$. Again, the
only thing we need is the residue class of just one factor of $n_m'$
modulo $n$ different from $1$ and the classes of the prime divisors
of $m+1$. Hence, if convenient, we can consider the set $S$ to be of
polynomial size such that $S\subset\{r\tm{n}\,:\,r\text{ prime },
r|n_m'\}$, $S\not\subseteq S_m\cup\{1\}$ where $S_m=\{r \tm{n}\,:\,r
\text{ prime }, r|(m+1)\}$. The following algorithm on the input of
an RSA modulus $n$ outputs a nontrivial factor of $n$.

\

\begin{alg}\label{det2}

\

\begin{enumerate}
\item  $m=2$
\item Send $(c,m)$  to $\cal O$.
\item If $S=\perp\Rightarrow m=m+1$ and go to (2). Else,
\item Take $r\in S$, $r\ne S_m\cup\{1\}$, and compute $d=(r-1,n)$.
\end{enumerate}
\end{alg}

\begin{thm}\label{gen}
Let $n=pq$ be an RSA modulus. If GRH is true then the Algorithm
\ref{det2} runs  in polynomial time and the number $d$ given by it
is a prime divisor of $n$.
\end{thm}
{\it Proof: } By Theorem \ref{sho} we can assume that $m$ is a
primitive root modulo $q$, at a polynomial time cost. Then
$m^{p-1}\not\equiv 1\tm q$, since $p<q$. Hence, in a similar way as
in the proof of Theorem \ref{male} we have to prove that a certain
prime factor $r$ of $n_m'$ belongs to a residue class modulo $n$ not
in $S_m\cup\{1\}$. We will use the following straightforward lemma.
\begin{lem}\label{prim}Let $n$ be an RSA modulus. For any integer $m$, such that $(m+1,n)=1$
we have $((m^n+1)/(m+1),m+1)=1$.
\end{lem}
{\it Proof:} Observe that if $r|(m+1)$, then
$$
(m^n+1)/(m+1)=\sum_{j=0}^{n-1}(-m)^j\equiv\sum_{j=0}^{n-1}1\tm
r\equiv n\tm r.
$$

\qed

Now, analogously to what we did in the proof of Theorem
\ref{male}, if $r|n_m'$ then $m^{2n}\equiv1\tm r$, and so
ord$_r(m)=2,p,q,2p,2q,pq$ or $2pq$ and clearly
ord$_r(m)\not=2$ for any $r$ a prime factor of $(m^{n}+1)/(m+1)$.
To see this use Lemma \ref{prim} and observe that if $r|(m-1)$
then $m^n+1\equiv 2\tm r$. Hence, as in the previous section, for
any $r|(m^n+1)/(m+1)$ then either $p|(r-1)$, $q|(r-1)$ or
$n|(r-1)$. If $r\equiv 1\tm n$ for any $r|(m^n+1)/(m+1)$ then
$m^{n-1}\equiv 1\tm n$, which is impossible for $m$ a primitive
root modulo $q$ since $m^{n-1}\equiv m^{p-1}\tm q$. The proof of
the theorem concludes as in Theorem \ref{male}.






\end{document}